\magnification=\magstep1
\input amstex
\documentstyle{amsppt}
\refstyle{C}
\PSAMSFonts
\pagewidth{13cm}
\pageheight{19cm}
\NoBlackBoxes
\newwrite\lab
\def\Label{\Labelunnecessary}   

\newcount\workdone
\workdone=2

\def \epip {(0.1)}
\def \tpiiip {Theorem 0.3}
\def \eipiip {(1.2)}
\def \eipviip {(1.7)}
\def \eipxvp {(1.15)}
\def \eipxvip {(1.16)}
\def \eiipiip {(2.2)}
\def \tiipiip {Theorem 2.2}
\def \liipiiip {Lemma 2.3}
\def \eiipivp {(2.4)}
\def \eiipviiip {(2.8)}
\def \eiipxvp {(2.15)}
\def \eiipxxip {(2.21)}
\def \eiipxxiip {(2.22)}
\def \eiipxxiiip {(2.23)}
\def \eiipxxvip {(2.26)}
\def \eiipxxviip {(2.27)}
\def \eiipxxxp {(2.30)}
\def \eiipxxxvip {(2.36)}
\def \liipivp {Lemma 2.4}
\def \eiipxliiip {(2.43)}
\def \eiiipivp {(3.4)}
\def \eiiipvp {(3.5)}
\def \eiiipviip {(3.7)}
\def \liiipip {Lemma 3.1}
\def \eiiipxp {(3.10)}
\def \eiiipxvp {(3.15)}
\def \eiiipxviiip {(3.18)}
\def \liiipiiip {Lemma 3.3}
\def \tiiipvp {Theorem 3.5}
\def \eiiipxxxip {(3.31)}
\def \eiiipxxxivp {(3.34)}
\def \eivpip {(4.1)}
\def \tivpiip {Theorem 4.2}
\def \livpiiip {Lemma 4.3}
\def \livpivp {Lemma 4.4}
\def \eivpxivp {(4.14)}
\def \livpviip {Lemma 4.7}
\def \eivpxxiiip {(4.23)}
\def \bk {1}
\def \bh {2}

\def \cgi {5}
\def \cgii {6}

\def \cgiv {8}
\def \giu {9}

\def \mm {14}
\def \tii {15}
\def \tiii {16}
\def \tw {17}


\ifnum \workdone=1 \immediate\closeout\lab\input \Label
\else\ifnum\workdone=0
\immediate\openout\lab=\Label\fi\fi


\define\apr#1#2#3#4{\csname\expandafter#1\romannumeral#2p%
           \romannumeral#3p\romannumeral#4\endcsname}
\define\Apr#1#2#3#4{\csname\string{#1\romannumeral#2p%
           \romannumeral#3p\romannumeral#4}\endcsname}


\define\wt#1#2#3{\ifnum \workdone=0 %
\immediate\write\lab{\csname def\endcsname\apr t{#1}{#2}{#3}{Theorem
\pc}}\else \fi}

\define\wl#1#2#3{\ifnum \workdone=0 
\immediate\write\lab{\csname def\endcsname\apr l{#1}{#2}{#3}{Lemma
\pc}}\else \fi}

\define\wD#1#2#3{\ifnum \workdone=0 
\immediate\write\lab{\csname def\endcsname\apr d{#1}{#2}{#3}{Definition
\pc}}\else \fi}

\define\wR#1#2#3{\ifnum \workdone=0 
\immediate\write\lab{\csname def\endcsname\apr r{#1}{#2}{#3}{Remark
\pc}}\else \fi}

\define\we#1#2#3{\ifnum \workdone=0 %
\immediate\write\lab{\csname def\endcsname\apr e{#1}{#2}{#3}{\tc}}\else \fi}

\define\wref#1{\ifnum \workdone=0 \ifx\\#1\\\empty \relax \else %
\immediate\write\lab{\csname def\endcsname\csname #1\endcsname
{\nr}}\fi\else \fi}


\define\gobble#1{}

\define\ci#1{\ifnum \workdone=1 #1\else
\ifnum \workdone=2 #1\else
{\bf?}\expandafter\gobble\string#1{\bf?} \fi\fi}

\define\cit#1#2#3#4{\ifnum \workdone=1 \apr #1#2#3#4 \else
\ifnum \workdone=2 \apr #1#2#3#4 \else
{\bf?}\expandafter\gobble\string\apr #1#2#3#4 {\bf?}\fi\fi}

\define\Cite#1{\ifnum \workdone=0 \ci#1\else \cite{#1}\fi}


\define\R#1{\ifnum \workdone=1 \csname #1\endcsname
\ifnum \workdone=2 \csname #1\endcsname  \else #1 \fi\fi}







\loadbold


\define\ndots{\raise 0.47ex \hbox {,}\hskip0.06em\cdots %
       \raise 0.47ex \hbox {,}\hskip0.06em} 


\newskip\Csmallskipamount                                                
\Csmallskipamount=\smallskipamount
\newskip\Cmedskipamount
\Cmedskipamount=\medskipamount
\newskip\Cbigskipamount
\Cbigskipamount=\bigskipamount

\define\Csmallskip{\vskip\Csmallskipamount} 
\define\Cmedskip{\vskip\Cmedskipamount}

\newdimen\spt
\spt=0.5pt

\define\vs#1{\vskip#1pt}
\define\hs#1{\hskip#1pt}
\define\hsc#1{\hskip#1\spt}


\define\bitem{\begingroup\parindent=32pt\Cmedskip\par\hangafter-1\textindent}
\define\enditem{\hfill\endgroup}

\define\citem{\Csmallskip\advance\itemno by
1\item"{(\romannumeral\the\itemno})"}


\define\graph{\operatorname{graph}}

\define\pa{\partial}
\define\pde#1#2{\dfrac {\partial#1}{\partial#2}}
\define\pd#1#2#3{\dfrac {\partial#1}{\partial#2^#3}}
\define\pdc#1#2#3{\dfrac {\partial#1}{\partial#2_#3}}
\define\Pd#1#2#3{\dfrac {{\partial\hskip0.15em}^2#1}{\partial {x^ #2}\,\partial
{x^#3}}}
\define\PD#1#2#3{\dfrac {{\partial\hskip0.15em}^2#1}{\partial
#2\,\partial#3}}
\define\ra{\rightarrow}

\define\sd{\vartriangle}
\define\sq#1{\sqrt{1+|D#1|^2}}
\define\su{\subset}
\define\Su{\Subset}
\define\sub{\operatorname{sub}}

\define\ti{\times }


\newcount\refno
\refno=0
\newcount\itemno
\itemno=0
\newcount\headno
\headno=-1
\newcount\equationno
\equationno=0
\newcount\proclaimno
\proclaimno=0



\define\hn{{\global\advance\headno by 1}\number\headno}
\define\en{{\global\advance\equationno by 1}\number\equationno}
\define\pn{{\global\advance\proclaimno by 1}\number\proclaimno}
\define\rn{{\global\advance\refno by 1}\number\refno}
\define\ct{\tag{\the\headno.\en}} 
\define\tc{(\the\headno.\the\equationno)}  
\global\define\cp{\the\headno.\pn} 
\define\pc{\the\headno.\the\proclaimno}
\define\nh{\the\headno} 
\define\nr{\the\refno} 

\topmatter
\title Closed hypersurfaces of prescribed mean curvature in locally conformally
flat Riemannian manifolds
\endtitle
\rightheadtext{ Closed hypersurfaces of prescribed mean curvature }

\author  Claus Gerhardt
\endauthor
\address {Ruprecht-Karls-Universit\"at \hfil}
\linebreak\indent {Institut f\"ur Angewandte Mathematik\hfil}
\linebreak\indent {Im Neuenheimer Feld 294 \hfil}
\linebreak\indent {D-69120 Heidelberg\hfil}
\linebreak\indent Germany
\endaddress
\email {gerhardt\@math.uni-heidelberg.de}
\endemail
\subjclass 35
\endsubjclass

\abstract {We prove the existence of smooth closed hypersurfaces of prescribed
mean curvature homeomorphic to
$S^n$ for small $n, n\le6$, provided there are barriers.}
\endabstract
\endtopmatter

%
\head \hn. Introduction
\endhead 
In a complete $(n+1)$-dimensional manifold $N$ we want to find closed
hypersurfaces $M$ of {\it prescribed mean curvature}. To be more precise, let
$\Omega $ be a connected open subset of $N$, $f\in C^{0,1}(\bar \Omega)$, then we
look for a closed hypersurface $M\su \Omega $ such that

$$ H_{|_{M}}=f(x)\qquad \forall\, x\in M,
\ct
$$
\we{0}{1}{0}%
where $H_{|M}$ is the mean curvature, i.e. the sum of the principal
curvatures.

The existence of a generalized solution $M=\partial E$, where $E$ is a Caccioppoli
set minimizing an appropriate functional is easily demonstrated if the boundary
of $\Omega $ is supposed to consist of two components acting as barriers. For
small
$n$, $n\le 6,$ the generalized solution is also a classical one, since it is smooth,
$M\in C^{2,\alpha }$, and hence a solution of \ci
\epip; but nothing is known about its topological type.

We shall prove that in the case when $n\le 6$ and $N$ is locally conformally flat, or
more precisely, when in $\Omega $ the metric is conformally flat, smooth
solutions homeomorphic to $S^n$ exist.

We make the following definition
\definition{Definition \cp} Let $M_1$, $M_2$ be closed hypersurfaces in $N$
homeomorphic to $S^n$ and of class
$C^{2,\alpha }$ which bound an open, connected, relatively compact subset
$\Omega $. $M_1$, $M_2$ are called {\it barriers} for
$(H,f)$ if
$$ H_{|_{M_1}}\le f
\ct
$$ and
$$ H_{|_{M_2}}\ge f
\ct
$$ Here, the mean curvature of $M_1$ is calculated with respect to the normal that
points outside of $\Omega $ while the mean curvature of $M_2$ is calculated with
respect to the normal that points inside of $\Omega $.
\enddefinition
\remark{Remark \cp}In view of the weak Harnack inequality the barriers do not
touch each other, unless both coincide and have prescribed mean curvature $f$. In
this case $\Omega $ would be empty.
\endremark We shall consider such  a region $\Omega $ bounded by barriers
$M_1$, $M_2$ for $(H,f)$, where $f\in C^{0,1}(\bar\Omega )$ is given, and assume
that $\Omega $ is conformally equivalent to an open, bounded set in $\bold
R^{n+1}$. Furthermore, we suppose, if we identify $\Omega $ with its image in
$\bold R^{n+1}$, that the barriers $M_1$, $M_2$ can be considered as graphs over
$S^n$, i.e., after fixing the origin and having introduced Euclidean polar
coordinates $(x^\alpha )_{0\le \alpha \le n}$, where
$x^0=r$ represents the radial distance, each $M_i$ can be written as a graph
$$ M_i=\graph u_{i{|_{S^n}}}=\{(x,x^0)\; :\; x^0=u_i(x),\,x\in S^n\},
\ct
$$ where we use slightly ambiguous notation.

The polar coordinates can also be considered to be a coordinate system in $N$
covering $\Omega $; the metric in $N$ then has the form
$$ d\bar s_N^{2}=e^{2\psi }d\bar s_{\bold R^{n+1}}^2=e^{2\psi }(dr^2+r^2\sigma
_{ij}dx^idx^j),
\ct
$$ where $(\sigma _{ij})$ is the standard metric on $S^n$.

Under these assumptions we shall prove
\proclaim{Theorem \cp} Let $\Omega $, $M_1$, $M_2$, and $f$ satisfy the
assumptions stated above, then the problem
$$ H_{|_M}=f
\ct
$$ has a solution $M\su \bar \Omega $ of class $C^{2,\alpha }$ homeomorphic to
$S^n$, if $n\le 6$.
\endproclaim
\wt{0}{3}{0}%
\remark{Remark \cp}Neither the function $f$ nor its derivatives are supposed to
satisfy any sign conditions. Even the assumption on the smoothness of $f$ can be
relaxed; if $f$ is only bounded, then a solution $M$ of class $H^{2,p}$ would exist
for any finite $p$.
\endremark

The problem of finding closed hypersurfaces of prescribed mean curvature has
been considered by Bakelman and Kantor
\Cite{\bk} and Treibergs and Wei \Cite{\tw} in the case when $N=\bold R^{n+1}$
assuming that $f$ is positive and satisfies
$$
\dfrac{\partial}{\partial r}(rf)\le 0
\ct
$$ where $r$ is the geodesic distance to some fixed origin. In \Cite{\cgi} we
proved the existence of a convex solution in space forms under the assumption
$$ -K_Nf\bar g_{\alpha \beta }+f_{\alpha \beta }\le 0\qquad \text{in } \Omega
\ct
$$ where in addition $f$ is supposed to be positive if $K_N>0$. In all cases the
existence of barriers is required.

The paper is organized as follows: In Section~1 we derive the basic equations for
hypersurfaces in conformally flat spaces, in Section~2 we solve auxiliary
problems, the solutions of which converge to the desired solution as is shown in
Sections~3 and 4.

\head \hn. Notations and preliminary results
\endhead\equationno=0\proclaimno=0  Let $N$ be a complete
$(n+1)$-dimensional manifold and $M$ a closed hypersurface. Geometric
quantities in $N$ will be denoted by $(\bar g_{\alpha\beta })$, $(\bar R_{\alpha
\beta \gamma \delta })$, etc., and those in $M$ by $(g_{ij})$, $(R_{ijkl})$, etc. Greek
indices range from $0$ to $n$ and Latin from $1$ to $n$; the summation convention
is always used. Generic coordinate systems in $N$ resp. $M$ will be denoted by
$(x^\alpha )$ resp. $(\xi ^i)$. Covariant differentiation will simply be indicated by
indices, only in case of possible ambiguity they will be preceded by a semicolon,
i.e. for a function $u$ on $N$, $(u_\alpha )$ will be the gradient and $(u_{\alpha
\beta })$ the Hessian, but, e.g. the covariant derivative of the curvature tensor
will be abbreviated by $\bar R_{\alpha \beta \gamma \delta ;\epsilon }$. We also
point out that
$$
\bar R_{\alpha \beta \gamma \delta ;i}=\bar R_{\alpha \beta \gamma \delta
;\epsilon }x_i^\epsilon 
\ct
$$ with obvious generalizations to other quantities.

In local coordinates $x^\alpha $ and $\xi ^i$ the geometric quantities of the
hypersurface $M$ are connected through the following equations

$$ x_{ij}^\alpha =-h_{ij}\nu ^\alpha 
\ct
$$
\we{1}{2}{0}%
 the so-called {\it Gau{\ss} formula}. Here, and also in the sequel, a
covariant derivative is always a {\it full\/} tensor, i.e.
$$ x_{ij}^\alpha =x_{,ij}^\alpha -\Gamma _{ij}^kx_k^\alpha +\bar \Gamma _{\beta
\gamma }^\alpha x_i^\beta x_j^\gamma .
\ct
$$ The comma indicates ordinary partial derivatives.

In this implicit definition {\ci \eipiip} the {\it second fundamental form\/}
$(h_{ij})$ is taken with respect to $-\nu $.

The second equation is the {\it Weingarten equation}
$$
\nu _i^\alpha =h_i^kx_k^\alpha 
\ct
$$

Finally, we have the {\it Codazzi equation\/}
$$ h_{ij;k}-h_{ik;j}=\bar R_{\alpha \beta \gamma \delta }\nu ^\alpha x_i^\beta
x_j^\gamma x_k^\delta  
\ct
$$  and the {\it Gau{\ss} equation}
$$ R_{ijkl}=h_{ik}h_{jl}-h_{il}h_{jk}+\bar R_{\alpha \beta \gamma \delta
}x_i^\alpha x_j^\beta x_k^\gamma x_l^\delta 
\ct
$$

Assume now, that the metric in $N$ is (locally) conformal to the metric in $\bold
R^{n+1}$
$$ d\bar s_N^2=e^{2\psi }d\bar s_{\bold R^{n+1}}^2\,,
\ct
$$\we{1}{7}{0}%
or more precisely, assume that {\ci \eipviip} is valid in the region
$\Omega $, where we shall consider $\Omega
$ to be a subset of $N$ as well as $\bold R^{n+1}$ without changing the notation.
The same convention applies to hypersurfaces
$M$ contained in $\Omega $ where we distinguish the geometric quantities of $M$
considered as a submanifold of $\bold R^{n+1}$ by using the notation $\hat h_{ij}$,
$\hat g_{ij}$, $\hat \nu ^\alpha $, etc. The connection with the corresponding
quantities in $N$ is given by
$$ g_{ij}=e^{2\psi }\hat g_{ij}\,,
\ct
$$
$$
\nu ^\alpha =e^{-\psi }\hat \nu ^\alpha,
\ct
$$ and
$$ h_{ij}e^{-\psi }=\hat h_{ij}+\psi _\alpha \hat \nu ^\alpha \hat g_{ij}\,.
\ct
$$ Thus, the mean curvatures of $M$ in $N$ resp. $\bold R^{n+1}$ are related
through
$$ He^\psi =\hat H+n\psi _\alpha \hat \nu ^\alpha .
\ct
$$

Assume now, that $M$ can be written as a graph over $S^n$, i.e after introducing
polar coordinates $(x^\alpha )$ in $\bold R^{n+1}$, where $x^0=r$,
$$M=\graph u_{|_{S^n}}=\big\{\,(r,x)\; :\; r=u(x),\; x\in S^n\,\big\}
\ct
$$ The metric in $\bold R^{n+1}$ is then expressed as
$$ d\bar s_{\bold R^{n+1}}^2=dr^2+r^2\sigma _{ij}dx^idx^j
\ct
$$ where $(\sigma _{ij})$ is the metric of $S^n$; the induced metric of $M$ is
$$
\hat g_{ij}=u_iu_j+u^2\sigma _{ij}=u^2(\varphi _i\varphi _j+\sigma _{ij}),
\ct
$$ where $\varphi =\log u$\hs1, and its inverse is given by
$$
\hat g^{ij}=u^{-2}\Big(\sigma ^{ij}-\dfrac{\varphi ^i}{v}\dfrac{\varphi ^j}{v}\Big),
\ct
$$
\we{1}{15}{0}%
where $v=\sqrt {1+\sigma ^{ij}\varphi _i\varphi _j}\equiv
\sqrt{1+|D\varphi |^2}$, $\sigma ^{ij}$ is the inverse of $\sigma _{ij}$\hs{0.5}, and
$\varphi ^i=\sigma ^{ij}\varphi _j$.

The product $u\hat H$ can be represented as
$$ u\hat H=-D_i(a^i(D\varphi ))+\frac{n}{v}
\ct
$$
\we{1}{16}{0}%
where $a^i(D\varphi )=v^{-1}\varphi ^i$ and the divergence is
calculated with respect to $\sigma _{ij}$.

\head \hn. Existence of solutions to an auxiliary problem
\endhead\equationno=0\proclaimno=0  As in \Cite{\cgii, Section 4} we first find
solutions to an auxiliary problem that will converge to the final solution.

\definition {Definition \cp}A hypersurface $M_0$ is called a supersolution for
$(H,f)$ if
$$ H_{|_{M_0}}\ge f.
\ct
$$
\enddefinition

We, furthermore, recall our assumption that the boundary of $\Omega $ consists
of barriers $M_i$, $i=1,2$\,, for $(H,f)$ that can be written as graphs over $S^n$,
$M_i=\graph u_i$.

We can now formulate the auxiliary problem.
\proclaim{Theorem \cp}Let $M_0=\graph u_0$, where $u_1\le u_0\le u_2$, and
$u_0$ is a supersolution for $(H,f)$ with $f\in C^{0,1}(\bar
\Omega )$, then the problem
$$
\left\{\aligned &H=f-\gamma e^{-\mu u}[u-u_0]\\ &u_1\le u\le u_0
\endaligned\right.
\ct
$$
\we{2}{2}{0}%
has a solution $u\in C^{2,\alpha }(S^n)$ provided $\gamma $ and $\mu
$ are sufficiently large, $\mu =\mu (\Omega ,f)$ and
$\gamma =\gamma (\mu ,\Omega ,f)$. Here, the reference that a term depends on
$\Omega $ should also indicate that geometrical quantities of the ambient space
and the barriers are involved.
\endproclaim
\wt{2}{2}{0}%

\subhead \nh.1. $\bold C^1$-estimates
\endsubhead\newline\newline Let $u\in C^{3,\alpha }(S^n)$ be a solution of {\rm\ci
\eiipiip}, where  we first assume a slightly higher degree of smoothness so that
the classical maximum principle can be applied to estimate $Du$, or equivalently,
the quantity $v=\sqrt{1+|D\varphi |^2}$ in {\rm\ci
\eipxvp}.

\proclaim{Lemma \cp}Let $u\in C^{3,\alpha }(S^n)$ be a solution of {\rm\ci
\eiipiip}, then $v=\sqrt{1+|D\varphi |^2}$ is estimated by
$$ v\le \text{\rm const}(|Du_0|,|Df|,\gamma ,\mu )
\ct
$$ provided $\gamma $ and $\mu $ are sufficiently large.
\endproclaim
\wl{2}{3}{0}%
\demo{Proof}We transfer the equation in {\rm\ci \eiipiip} into $\bold R^{n+1}$, i.e.
we consider
$$
\hat H=e^\psi H-n\psi _\alpha \hat \nu =e^\psi \bigl\{f-\gamma e^{-\mu
u}[u-u_0]\bigr\}-n\psi _\alpha \hat \nu ^\alpha.
\ct
$$
\we{2}{4}{0}%
Let $x$ be the embedding vector for the hypersurface, then we
define
$$
\chi =\langle x,\hat \nu \rangle^{-1}=u^{-1}\langle Dr,\hat \nu \rangle^{-1}=u^{-1}v
\ct
$$ and we shall prove a priori estimates for $\chi $.

We choose local coordinates and compute the first and the second covariant
derivatives of $\chi $
$$
\chi _{\tsize_i}=-\langle x,\hat \nu \rangle^{-2}\langle x,\hat \nu _i\rangle=-\chi
^2\hat h_i^k\langle x_k,x\rangle,
\ct
$$
$$
\chi _{\tsize_{ij}}=2\chi ^{-1}\chi _{\tsize_i}\chi _{\tsize_j}-\chi ^2\hat
h_{i;j}^k\langle x_k,x\rangle+\chi \hat h_i^k\hat h_{kj}-\chi ^2\hat h_{ij}.
\ct
$$ Hence, we conclude in view of the Codazzi equations
$$ -\Delta \chi =-\hat g^{ij}\chi _{\tsize_{ij}}=-2\chi ^{-1}\|D\chi \|^2-|\hat A|^2\chi
+\hat H\chi ^2+\hat H_ku^ku\chi ^2.
\ct
$$
\we{2}{8}{0}%
Here, all indices are raised with respect to the induced metric, and
we used the abbreviations
$$
\|D\chi \|^2=\hat g^{ij}\chi _{\tsize_i}\chi _{\tsize_j}\qquad \text{and}\qquad |\hat
A|^2=\hat h_i^k\hat h_k^i.
\ct
$$

The crucial terms are those which are quadratic in $\chi $, they have to add up to
something negative, if $\chi $ is large.

To compute $\hat H_k$ it is convenient to introduce polar coordinates $(x^\alpha
)$ in $\bold R^{n+1}$ and to decompose
$D_\alpha f\/$ into its radial part and into the tangential components with respect
to $S^n$
$$
\dot f=\pde fr\qquad \text{and}\qquad f_{,i}=\pd fxi \,.
\ct
$$

We then obtain from {\rm\ci \eiipivp}
$$
\align
\hat H _k=&\medspace e^\psi \big\{f-\gamma e^{-\mu u}[u-u_0]\big\}\psi _\alpha
x_k^\alpha -n\psi _{\alpha \beta }\hat \nu ^\alpha x_k^\beta -n\psi _\alpha
\hat h_k^ix_i^\alpha
\ct\\ &+e^\psi \bigl\{\dot fu_k+f_{,k}+\gamma \mu e^{-\mu u}[u-u_0]u_k-\gamma
e^{-\mu u}[u_k-u_{0,k}]\bigr\}
\endalign
$$ Using the relations
$$ u^k=u^{-2}\sigma ^{ik}u_iv^{-2},\quad \|Du\|^2=\frac {|D\varphi |^2}{v^2}
\ct
$$ and
$$\chi _{\tsize_i}=-\chi ^2\hat h_i^ku_ku
\ct
$$ we conclude
$$\multline
\hat H_ku^ku\le \medspace c|f|+c|Df|-c\gamma e^{-\mu u}[u-u_0]\big[\|Du\|^2+\chi
^{-2}\big]+c\chi ^{-1}+\chi ^{-3}\|D\chi \|^2\\ +\gamma \mu e^{-\mu
u}[u-u_0]\|Du\|^2ue^\psi
 -\frac{1}{2}\gamma e^{-\mu u}\bigl[\|Du\|^2-u^{-2}v^{-2}\sigma
^{ij}D_iu_0D_ju_0\bigr]ue^\psi 
\endmultline\tag{\nh.\en}
$$ where $c=c(\Omega )$. Therefore, the righthand-side in {\rm\ci \eiipviiip} can
be estimated from above by
$$
\multline c\bigl[|f|+|Df|+1\bigr]\chi ^2+(\mu -c)\|Du\|^2\gamma e^{-\mu
u}[u-u_0]ue^\psi \chi ^2\\ -\frac{1}{2}\gamma e^{-\mu
u}\bigl[\|Du\|^2-u^{-2}v^{-2}\sigma ^{ij}D_iu_0D_ju_0\bigr]ue^\psi\chi ^2
\endmultline\tag{\nh.\en}
$$
\we{2}{15}{0}%
at points where $\|Du\|^2\ge \frac{1}{2}$. We now choose $\mu $
larger than $2c$, so that the second term in {\rm\ci
\eiipxvp} is non-positive. Then, we choose $\gamma $ such that
$$ c\bigl[|f|+|Df|+1\bigr]\le \frac{1}{4}\gamma e^{-\mu u}ue^\psi 
\ct
$$ and deduce from the maximum principle that
$$ v\le c(\Omega ,Du_0)
\ct
$$ where the constant is determined by the relation
$$ |D\varphi |^2-u^{-2}\sigma ^{ij}D_iu_0D_ju_0\le \frac{1}{2}v^2
\ct
$$
\enddemo

\subhead
\nh.2. Existence of solutions
\endsubhead\newline\newline We are still looking for hypersurfaces in $\bold
R^{n+1}$, i.e. we want to solve the equation {\rm\ci{\eiipivp}} with the
side-conditions
$$ u_1\le u\le u_0,
\ct
$$ or equivalently, we can solve
$$ -D_i(a^i(D\varphi ))+ \frac{n}{v}=u\hat H=ue^\psi \bigl\{f-\gamma e^{-\mu
u}[u-u_0]\bigr\}-un\psi _\alpha \hat \nu ^\alpha ,
\ct
$$ where $\varphi =\log u$\hs1, cf. equation {\rm\ci \eipxvip}.

Let us denote the lower order terms in the preceding equation by $a(x,\varphi
,D\varphi )$, and let $\varphi _1=\log u_1$\hs1, $\varphi _0=\log u_0$\hs1, then
we have to solve
$$
\left\{
\aligned -&D_i(a^i(D\varphi ))+a(x,\varphi ,D\varphi )=0\\ &\varphi _1\le \varphi
\le \varphi _0
\endaligned
\right.
\ct
$$
\we{2}{21}{0}%
Here, $a^i$ is a strictly monotone vectorfield, the lower order
term and its derivatives are bounded in $\varphi _1\le \varphi
\le \varphi _0$, i.e.
$$ |a|+\Big|\pd axi\Big|+\Big |\pde a\varphi \Big|+\Big|\pd api\Big|\le \text{const},
\ct
$$
\we{2}{22}{0}%
and moreover
$$
\pde a\varphi \ge \epsilon_{\tsize_0}>0\qquad \text{in}\qquad \varphi _1\le
\varphi \le \varphi _0,
\ct
$$
\we{2}{23}{0}%
due to our choice of $\mu $ and $\gamma $, where we increase $\mu
$ and $\gamma $ a bit in view of the presence of the additional factor $e^\varphi $.

To solve {\rm\ci \eiipxxip} we first assume that $f$ is of class $C^{1,\alpha }$ and
$u_0$ of class $C^{3,\alpha }$. Extend
$f_0=H_{|_{M_0}}$, where $M_0=\graph u_0$, to $\bar \Omega $ by setting
$$ f_0(x,r)=f_0(x)\; ,\qquad x\in S^n
\ct
$$ and consider the convex combination
$$ f_t=tf+(1-t)f_0\; ,\qquad 0\le t\le 1.
\ct
$$

Then, we look at the problems
$$
\left\{
\aligned -&D_i(a^i(D\varphi _t ))+a(x,\varphi _t,D\varphi _t)=0\\ &\varphi _1\le
\varphi _t\le \varphi _0
\endaligned
\right.
\ct
$$
\we{2}{26}{0}%
where $f$ is replaced by $f_t$, and where we have a slight
ambiguity in the notation for $t=1$. The lower order term also depends explicitly
on $t$, but since the estimates {\rm\ci \eiipxxiip} and {\rm\ci \eiipxxiiip} are
independent of $t$, if we choose $\gamma $ sufficiently large---at the moment
$\gamma $ also depends on $f_0$---, we do not indicate it specifically.

We shall use the continuity method to prove that {\rm\ci \eiipxxvip} has a
solution for all $0\le t\le 1$. Let us treat {\rm\ci
\eiipxxvip} as a variational inequality
$$
\left\{
\aligned
\langle-&D_i(a^i(D\varphi _t))+a(x,\varphi _t,D\varphi _t),\eta -\varphi
_t\rangle\ge 0\qquad \forall\, \eta \in K,\\ K&=\big\{\,\eta \in C^{0,1}(S^n)\; :\;
\varphi _1\le \eta \le \varphi _0\,\big\}.
\endaligned
\right.
\ct
$$
\we{2}{27}{0}%
It can easily be shown that the obstacles $\varphi _1$, $\varphi _0$
act as barriers, i.e. they are sub- resp. supersolutions for any value of $t$, $0\le
t\le 1$, so that any solution of the variational inequality is actually a solution of
the corresponding equation. However, our proof of the solvability of {\rm\ci
\eiipxxviip} is valid for arbitrary $C^{2,\alpha }$ obstacles.

Define $\Lambda $ through
$$
\Lambda =\big\{\,t\in [0,1]\; :\; \text{{\rm\ci \eiipxxviip} has a solution }  \varphi
_t\,\big\}.
\ct
$$ Then, we conclude

\bitem{(i)} $\Lambda \ne \emptyset$, since $0\in \Lambda $.\enditem
\bitem{(ii)} $\Lambda $ {\it is closed.}  It is well known that any solution of the
variational inequality is of class
$H^{2,p}(S^n)$ for any finite $p$. Therefore, the solution does not touch the
obstacles at points, where the gradient is larger than the gradients of the
obstacles, and it is a solution of the equation there. At those points the solution 
is also of class
$C^{3,\alpha }$ because $f$, $f_0$ are of class $C^{1,\alpha }$ and {\rm\ci \liipiiip} is
thus applicable. We conclude further that uniform $H^{2,p}$-estimates are valid
and hence uniform $C^{1,\alpha }$-estimates, which proves the closedness of
$\Lambda $.
\enditem

\bitem{(iii)}$\Lambda $ {\it is open.} To prove the openness we argue as in a
former paper {\rm\Cite{\cgiv}}.  Let $t_0\in
\Lambda
$ and let $\varphi _{t_0}$ be the corresponding solution of {\rm\ci \eiipxxviip}
with
$$ |D\varphi _{t_0}|\le c_0
\ct
$$
Let $\tilde a^i=\tilde a^i(p)$ be a uniformly monotone vectorfield that
coincides with $a^i$ for $|p|\le c_0+1$. The existence of such a vectorfield has
been shown in {\rm\Cite{\cgiv, Appendix II}}. Then, the corresponding variational
inequality\enditem
$$
\left\{
\aligned
\langle-&D_i(\tilde a^i(D\tilde \varphi _t))+a(x,\tilde \varphi _t,D\tilde \varphi
_t), \eta -\tilde \varphi _t\rangle\ge 0\qquad
\forall\, \eta \in \widetilde K\\
\widetilde K&=\big\{\,\eta \in H^{1,2}(S^n)\; :\; \varphi _1\le \eta \le \varphi
_0\,\big\}
\endaligned
\right.
\ct
$$
\we{2}{30}{0}%
has a solution $\tilde \varphi _t\in H^{2,p}(S^n)$ for any $t$, $0\le
t\le 1$, since the differential operator is uniformly elliptic in
$\widetilde K$, and there exists $\lambda >0$ such that the operator
$$A\varphi =-D_i(\tilde a^i(D\varphi ))+a(x,\varphi ,D\varphi )+\lambda \varphi 
\ct
$$ is uniformly monotone, i.e. there exists $\epsilon _{\tsize_0}>0$ such that
$$
\epsilon _{\tsize_0}\|\varphi -\eta \|_{1,2}^2\le \langle A\varphi -A\eta , \varphi
-\eta \rangle\qquad \forall\, \varphi , \eta \in
\widetilde K,
\ct
$$ where the norm on the left-hand side is the norm in $H^{1,2}(S^n)$. The operator
$$ -D_i(\tilde a^i(D\varphi ))+a(x,\varphi ,D\varphi )
\ct
$$ is therefore {\it pseudomonotone} and {\it coercive} in $\widetilde K$, and the
existence of solutions for the problem {\rm\ci
\eiipxxxp} follows from the general theory for solutions of variational
inequalities, cf. {\rm\Cite{\bh}}.

As we shall show below, the solutions of {\rm\ci \eiipxxxp} are unique for each
$t$, hence, they depend continuously on $t$ in the
$C^{1,\alpha }$-norm, and we conclude that for small $\epsilon >0$
$$ |D\tilde \varphi _t|\le c_0+\frac{1}{2}\qquad \forall\, t\in B_\epsilon (t_0)
\ct
$$ since $\tilde \varphi _{t_0}=\varphi _{t_0}$, and we deduce further, that these
$\tilde \varphi _t$ are also solutions of {\rm\ci
\eiipxxviip}, i.e. $\Lambda $ is open.

We have thus proved that $\Lambda $ coincides with the whole interval, so we
have especially proved the existence of a solution for the crucial value $t=1$. In
this case, $f_t=f$ and the obstacles are barriers, so that we deduce with the help
of the weak Harnack inequality, that the solution of the variational inequality is
actually a solution of the equation. For details see the uniqueness proof in {\rm\ci
\liipivp} below.

Let us point out, that at the moment the parameter $\gamma $ also depends on
$f_0$, and hence on the second derivatives of
$u_0$. However, $\gamma $ should only depend on $|u_0|$ and on the other
quantities mentioned in {\rm\ci \liipiiip}. To achieve this result, we more or less
repeat the argument just given in the first part of the existence proof.

Let $\gamma _0$ be a constant such that the gradient estimate in {\rm\ci \liipiiip}
and the relation {\rm\ci \eiipxxiiip} are valid for $\gamma \ge \gamma_0$. Let
$\bar \gamma \ge \gamma_0$ be arbitrary and define $\Lambda $ through
$$
\Lambda =\big\{\,\gamma \ge \bar \gamma \; :\; \text{{\rm\ci \eiipxxip} has a
solution}\,\big\}.
\ct
$$
$\Lambda $ is not empty as we have just proved. Let $\gamma ^\star=\inf \gamma
$. By repeating the arguments we used to prove that the variational inequality has
a solution, we conclude that $\gamma ^\star=\bar \gamma $, i.e. the existence of a
solution to the auxiliary problem is guaranteed provided $\gamma \ge \gamma
_0$.

Before we prove the uniqueness of the solution to the variational inequality, let
us remove the additional assumptions regarding the smoothness of $f$ and $u_0$.
We assumed in the proof $f\in C^{1,\alpha }(S^n)$ and $u_0\in C^{3,\alpha }(S^n)$, so
that the solutions $\varphi _t^{}$ of {\rm\ci \eiipxxviip} are of class $C^{3,\alpha }$
at points where they do not touch the obstacles in order to apply the classical
maximum principle to estimate the $C^1$-norm. But the gradient estimate for the
final solution, when $t=1$, only depends on the $C^1$-norms of $f$ and $u_0$,
hence, we obtain solutions under the weaker assumptions by approximation.

To complete the proof of the theorem, let us now show
\proclaim{Lemma \cp}Let $\varphi \in K$ be a solution of the variational inequality
$$
\left\{
\aligned
\langle-&D_i(a^i(D\varphi ))+a(x,\varphi ,D\varphi ),\eta -\varphi \rangle\ge
0\qquad \forall\, \eta \in K,\\ K&=\big\{\,\eta \in C^{0,1}(S^n)\; :\; \varphi _1\le
\eta \le \varphi _0\,\big\},
\endaligned
\right.
\ct
$$
\we{2}{36}{0}%
then $\varphi $ is uniquely determined, where we assume that the
condition {\rm\ci \eiipxxiiip} is valid and $a^i$, $a$ are of class $C^1$ in their
arguments.
\endproclaim
\wl{2}{4}{0}%
\demo {Proof}Let $\varphi $, $\tilde \varphi $ be two solutions of {\rm\ci
\eiipxxxvip}, then we have to show $\varphi =\tilde
\varphi $, or by symmetry, $\varphi \ge \tilde \varphi $.

We know that $\varphi $, $\tilde \varphi $ are of class $H^{2,p}(S^n)$ for any finite
$p$. Suppose
$$ G=\big\{\,x\; :\; \varphi <\tilde \varphi \,\big\}\ne \emptyset,
\ct
$$ then $G$ is open and in $G$ we have
$$ -D_i(a^i(D\varphi ))+a(x,\varphi ,D\varphi )\ge 0\qquad \text{because}\qquad
\varphi <\varphi _2
\ct
$$ and
$$ -D_i(a^i(D\tilde \varphi ))+a(x,\tilde \varphi ,D\tilde \varphi )\le 0\qquad
\text{because}\qquad \tilde \varphi >\varphi _1.
\ct
$$ Hence, we infer
$$ -D_i(a^i(D\varphi ))+D_i(a^i( D\tilde \varphi )) +a(x,\varphi ,D\varphi
)-a(x,\tilde \varphi , D\tilde \varphi )\ge 0\,,
\ct
$$ or, by setting $\varphi _t^{}=t\varphi +(1-t)\tilde \varphi $, $0\le t\le 1$, and
using the main theorem of calculus
$$ -D_i(\hat a^{ij}D_j(\varphi -\tilde \varphi ))+\pde {\hat a}{\varphi }(\varphi
-\tilde \varphi )+\pdc {\hat a}pi D_i(\varphi -\tilde
\varphi )\ge 0\,,
\ct
$$ where
$$
\aligned
\hat a^{ij}&=\int_0^1a^{ij}(D\varphi _t)\,,\\
\pde{\hat a}\varphi &=\int_0^1 \pde a\varphi (x,\varphi _t,D\varphi _t)\,,\\
\pdc {\hat a}pi&=\int_0^1 \pdc api(x,\varphi _t,D\varphi _t)\,,
\endaligned
\ct
$$ and we conclude that
$$ -D_i(\hat a^{ij}D_j(\varphi -\tilde \varphi ))+\pdc {\hat a}pi D_i(\varphi
-\tilde\varphi )\ge -\pde {\hat a}{\varphi }(\varphi -\tilde\varphi )
\ct
$$
\we{2}{43}{0}%
in $G$.

Now, by assumption
$$ m_0=\inf (\varphi -\tilde \varphi )<0\,,
\ct
$$ and the infimum is attained at a point $x_0\in G$. Define $\eta =\varphi -\tilde
\varphi -m_0$, then $\eta \ge 0$\hs{0.5}, and {\rm\ci \eiipxliiip} yields
$$ -D_i(\hat a^{ij}D_j\eta )+\pdc{\hat a}piD_i\eta >0
\ct
$$contradicting the weak Harnack inequality that would demand $\eta \equiv 0$.
Thus, we deduce $m_0\ge 0$\hsc1, and the uniqueness is proved.
\enddemo

\head \hn. Almost minimal solutions
\endhead\equationno=0\proclaimno=0 

We now apply the existence result of {\rm\ci \tiipiip} successively. Let $u_2$ be
{\it the} upper barrier; then, if $u_{k-1}$ is already defined for $k\ge 3$, let
$u_k\in C^{2,\alpha }(S^n)$ be the solution of
$$
\left\{
\aligned &H=f-\gamma e^{-\mu u_k}[u_k-u_{k-1}]\\ &u_1\le u_k\le u_{k-1}
\endaligned
\right.
\ct
$$ The solutions $(u_k)$ form a bounded monotone decreasing sequence, which
converges pointwise to a function $u$. The mean curvatures of the graphs
converge pointwise to $f(x,u)$, since $\gamma $ and $\mu $ are fixed; hence,
$\graph u$ would be a solution of our problem, if the $u_k$'s would satisfy uniform
$C^1$-estimates. But unfortunately, we cannot prove this, it might even be false.
Gradient estimates for graphs depend on the Lipschitz constant of the mean
curvature, i.e. $|Du_k|$ depends on $|Du_{k-1}|$.

However, the regularity results of De Giorgi, Massari, and Tamanini for almost
minimal hypersurfaces imply uniform
$C^{1,\frac{1}{2}}$-estimates for the {\it hypersurfaces} provided the
hypersurfaces are almost minimal, their mean curvatures uniformly bounded, and
$n\le 6$, cf. {\rm\Cite{\tii, \ci\tiii}}.

To apply these results, we shall prove that the hypersurfaces $M_k=\graph \log
u_k$ are almost minimal in the {\it metric product} $S^n\times \bold R$.

We adopt the view point and the notations from Section 2, i.e. we consider the
hypersurfaces as submanifolds of $\bold R^{n+1}$ and look at their diffeomorphic
images in $S^n\times \bold R$ under the diffeomorphism $\Phi (x,r)=(x,\log r)$.
Then, each $\varphi _k=\log u_k$ satisfies the equation in {\rm\ci \eiipxxip} on
$S^n$. For notational reasons we drop the index $k$, having in mind that it is fixed.
Furthermore, we consider the lower order term $a(x,\varphi ,D\varphi )$ to
depend only on $x$ without changing the symbol, i.e. $\varphi $ is a solution of 
$$ -D_i(a^i(D\varphi ))+a(x)=0\qquad \text{in }  S^n,
\ct
$$ where $a(x)$ is uniformly bounded.

Instead of $a(x)$ let us consider the modified lower order term
$$ a_{\epsilon _0}(x,t)=a(x)+\epsilon _{\tsize_0}(t-\varphi (x))
\ct
$$ with $\epsilon _{\tsize_0}>0$. Then, $a_{\epsilon _0}(x,\varphi )=a(x)$ and
therefore, we have
$$ -D_i(a^i(D\varphi ))+a_{\epsilon _0}(x,\varphi )=0
\ct
$$
\we{3}{4}{0}%
with
$$
\pde{a_{\epsilon _0}}\varphi =\epsilon _{\tsize_0}>0.
\ct
$$
\we{3}{5}{0}%
We shall prove that the boundary of the subgraph
$$ E=\sub \varphi =\big\{\,(x,t)\; :\; t<\varphi (x),\;  x\in S^n\,\big\}
\ct
$$ is almost minimal in the metric product $S^n\times \bold R$, or more precisely,
that it solves the variational problem of minimizing the so-called {\it perimeter}
plus an additional prescribed mean curvature term in $S^n\times \bold R$.

We first show that $E$ is minimal compared with other subgraphs.
\proclaim{Lemma \cp}The solution $\varphi $ of {\rm\ci \eiiipivp} is also a
solution of the variational problem
$$ J(\eta )=\int_{S^n}\sqrt{1+|D\eta |^2} +\int_{S^n}\int_0^\eta a_{\epsilon
_0}(x,t)\ra \min\qquad \forall\,\eta \in BV(S^n).
\ct
$$
\we{3}{7}{0}%
\endproclaim
\wl{3}{1}{0}%

$BV(S^n)$ is the space of functions of bounded variation, i.e. functions the
derivatives of which are bounded measures. For $\eta
\in BV(S^n)$ the area term in {\rm\ci \eiiipviip} is defined by
$$
\int_{S^n}\sqrt{1+|D\eta |^2}=\sup \Big\{\,\int_{S^n}(\gamma ^0+\eta D_i\gamma
^i)\; :\; \gamma ^\alpha \in C^\infty (S^n),\; |\gamma ^0|^2+\sigma _{ij}\gamma
^i\gamma ^j\le  1\Big\}.
\ct
$$ It coincides with the usual definition if $\eta $ is Lipschitz continuous.
\demo{Proof {\rm  of} {\rm\ci \liiipip}}The functional in {\rm\ci \eiiipviip} consists
of the standard area for graphs plus a mean curvature term; the corresponding
Euler-Lagrange equation is exactly the equation in {\rm\ci \eiiipivp}, thus, it is
not surprising that $\varphi $ should also solve {\rm\ci \eiiipviip}, since we know
from {\rm\ci \liipivp}, that the solutions of {\rm\ci
\eiiipivp} are uniquely determined.

Let $c_0$ be an arbitrary constant such that
$$ |D\varphi |<c_0.
\ct
$$ Then, solve the variational problem
$$
\left\{
\aligned &J(\eta )\ra \min\qquad \forall\, \eta \in K,\\ &K=\big\{\, \eta \in
C^{0,1}(S^n)\; :\; |D\eta |\le c_0\,\big\}.
\endaligned
\right.
\ct
$$
\we{3}{10}{0}%
Let $\tilde \varphi $ be a solution of {\rm\ci \eiiipxp}, the
existence of which can easily be proved in view of {\rm\ci \eiiipvp}, then $\tilde
\varphi $ solves the variational inequality
$$
\langle-D_i(a^i(D\tilde \varphi ))+a(x,\tilde \varphi ),\eta -\tilde \varphi
\rangle\ge 0\qquad \forall\, \eta \in K.
\ct
$$ On the other hand, since $\varphi $ is a solution of the equation and belongs to
$K$, we deduce from the strict monotonicity of the operator that $\varphi =\tilde
\varphi $.

Thus, $\varphi $ is a solution of the unconstrained variational problem
$$ J(\eta )\ra \min\qquad \forall\, \eta \in C^{0,1}(S^n),
\ct
$$ since $c_0$ is arbitrary, and by approximation we conclude, that $\varphi $ also
minimizes the functional in $BV(S^n)$.
\enddemo

Let $\widetilde N=S^n\times \bold R$ be the metric product of $S^n$ and $\bold
R$, so that the metric in $\widetilde N$ is given by
$$ ds_{\tilde N}^2=dt^2+\sigma _{ij}dx^idx^j.
\ct
$$ We also use $x^0$ instead of $t$ when appropriate.

The {\it perimeter} of a measurable set $E\su \widetilde N$ with respect to an
open set $\widetilde \Omega \su \widetilde N$ is defined by
$$
\int_{\tilde \Omega }|D\chi _{\tsize_E}|=\sup \Big\{\int_{\tilde \Omega }\chi
_{\tsize_E}D_\alpha \gamma ^\alpha \; :\;
\gamma ^\alpha \in C_c^\infty(\widetilde \Omega ),\;  |\gamma ^0|^2+\sigma
_{ij}\gamma ^i\gamma ^j\le 1\,\Big\},
\ct
$$ i.e. $E$ has finite perimeter in $\widetilde \Omega $ iff $\chi _{\tsize_E}$
belongs to $BV(\widetilde \Omega )$. Sets of finite perimeter are also called {\it
Caccioppoli} sets. It is well known that the perimeter of subgraphs is equal to the
area of the boundary.
\proclaim{Lemma \cp}Let $\Omega \su S^n$ be open, $\eta \in BV(\Omega )$ and
$E=\sub \eta _{|_\Omega }$, then
$$
\int_\Omega \sqrt{1+|D\eta |^2}=\int_{\Omega \times \bold R}|D\chi _{\tsize_E}|
\ct
$$
\we{3}{15}{0}%
\endproclaim

The proof is the same as in the case when $E$ is a subgraph in $\bold R^{n+1}$, cf.
{\rm\Cite{\giu, Theorem 14.6}}; moreover, we only need the relation when $\eta $
is of class $C^1$ and then {\rm\ci \eiiipxvp} follows immediately from the
divergence theorem.

The demonstration of the next lemma is also identical to the proof of its
Euclidean counterpart which is due to Miranda, cf. {\rm\Cite{\mm}}, but for the
convenience of the reader we shall repeat a version of the proof that can be found
in {\rm\Cite{\giu, Lemma 14.7}}.
\proclaim{Lemma \cp}Let $\Omega \su S^n$ be open, $F\su \Omega \times \bold
R$ be measurable such that
$$
\Omega \times (-\infty,-T)\su F\su \Omega \times (-\infty,T).
\ct
$$ For $x\in \Omega $ define
$$
\psi (x)=\lim_{k\ra \infty}\Big\{\int_{-k}^{k}\chi _{\tsize_F}(x,t)-k\Big\}.
\ct
$$ Then, there holds
$$
\int_\Omega \sqrt{1+|D\psi |^2}\le \int_{\Omega \times \bold R}|D\chi _{\tsize_F}|.
\ct
$$
\we{3}{18}{0}%
\endproclaim
\wl{3}{3}{0}%
\demo{Proof} We note that $\partial F\cap \Omega \times \bold R\su \Omega
\times (-T,T)$. Let
$$
\psi _k=\int_{-k}^k\chi _{\tsize_F}(x,t)-k,
\ct
$$ then $\psi _k$ is stationary for $k\ge T$, for 
$$
\int_{-k}^k\chi _{\tsize_F}(x,t)=\int_{-T}^T\chi
_{\tsize_F}(x,t)+\int_{-k}^{-T}1=\int_{-T}^T\chi _{\tsize_F}(x,t)-T+k,
\ct
$$ i.e.
$$
\psi =\int_{-T}^T\chi _{\tsize_F}(x,t)-T
\ct
$$ and $-T\le \psi \le T$.

Consider now $\gamma ^\alpha \in C_c^\infty(\Omega )$, $0\le \alpha \le n$,
satisfying $|\gamma ^0|^2+\sigma _{ij}\gamma ^i\gamma ^j\le 1$, and a smooth
real function $\eta $ such that $0\le \eta \le 1$ and
$$
\left\{
\aligned &\eta (t)=0\; ,\qquad |t|\ge T+1,\\ &\eta (t)=1\; ,\qquad |t|\le T.
\endaligned
\right.
\ct
$$ Then, we derive
$$
\int _{-\infty}^\infty\dot \eta \chi _{\tsize_F}=\int_{-\infty}^{-T}\dot \eta =\eta
(-T)=1,
\ct
$$ and
$$
\int _{-\infty}^\infty\eta \chi _{\tsize_F}=\int_{-T-1}^{-T} \eta+\int_{-T}^T\chi
_{\tsize_F}=\psi +T+\int_{-T-1}^{-T}\eta
\equiv \psi +c
\ct
$$ with $c=\text{const}$, from which we infer
$$
\align
\int_{\Omega \ti\bold R}|D\chi _{\tsize_F}|&\ge \int_{\Omega \ti\bold R}\chi
_{\tsize_F}D_\alpha (\eta \gamma ^\alpha )\tag{\nh.\en}\\ &=\int_\Omega
\int_{-\infty}^\infty \chi _{\tsize_F}\dot \eta \gamma ^0+\int_\Omega
\int_{-\infty}^\infty \chi _{\tsize_F}\eta D_i\gamma ^i\\ &=\int_\Omega \gamma
^0 +\int_\Omega (\psi +c)D_i\gamma ^i
\endalign
$$ and hence
$$
\int_{\Omega \ti \bold R}|D\chi _{\tsize_F}|\ge \int_\Omega \sq \psi .
\ct
$$
\enddemo

Let us return to the solution $\varphi $ of {\rm\ci \eiiipivp} that also solves the
variational problem {\rm\ci \eiiipviip}. We are going to prove that $E=\sub
\varphi $ locally minimizes the functional
$$
\Cal F(F,\widetilde \Omega )=\int_{\tilde \Omega }|D\chi _{\tsize_F}|+\int_{\tilde
\Omega }\chi _{\tsize_F}a_{\epsilon _0}
\ct
$$ which is defined for any $\widetilde \Omega \Su \widetilde N$ and any
Caccioppoli set $F\su \widetilde N$.
\definition{Definition \cp}A Caccioppoli set $E\su \widetilde N$ is said to be a
local minimizer for the functional $\Cal F$ if for any open $\widetilde \Omega \Su
\widetilde N$ and any Caccioppoli set $F$ with $F\vartriangle E\Su \widetilde
\Omega $ we have 
$$
\Cal F(E,\widetilde \Omega )\le \Cal F(F,\widetilde \Omega ).
\ct
$$
\enddefinition

We are now ready to prove
\proclaim{Theorem \cp}Let $\varphi $ be the solution of the variational problem
{\rm\ci \eiiipviip} and $E=\sub \varphi $, then
$E$ is a local minimizer for the functional $\Cal F$.
\endproclaim
\wt{3}{5}{0}%
\demo{Proof}Let $\widetilde \Omega \Su \widetilde N$ and let $F$ be a
Caccioppoli set with $F\sd E\Su \widetilde \Omega $, then $F$ satisfies the
conditions of {\rm\ci \liiipiiip} since $\varphi $ is bounded---actually any
$BV(S^n)$ solution of {\rm\ci
\eiiipviip} is bounded as it is well known---, and where we choose $\Omega =S^n$.

Define $\psi $ as in {\rm\ci \liiipiiip} and set $F^\star=\sub \psi $; then, we deduce
from {\rm\ci \eiiipviip}, {\rm\ci
\eiiipxvp}, and {\rm\ci \eiiipxviiip}
$$
\int_{\tilde \Omega }|D\chi _{\tsize_E}|\le \int_{\tilde \Omega }|D\chi
_{\tsize_F}|+\int_{S^n}\int_\varphi ^\psi a_{\epsilon _0}(x,t).
\ct
$$ We now observe that for arbitrary but fixed $x\in S^n$
$$
\int_\varphi ^\psi a_{\epsilon _0}(x,t)=\int_{-k}^ka_{\epsilon _0}(x,t)[\chi
_{\tsize_{F^\star}}-\chi _{\tsize_E}],
\ct
$$ where $k\ge |\varphi |+|\psi |$, and we claim furthermore, that
$$
\int_{-k}^ka_{\epsilon _0}(x,t)[\chi _{\tsize_{F^\star}}-\chi _{\tsize_F}]\le 0,
\ct
$$
\we{3}{31}{0}%
since $a_{\epsilon _0}(x,\hs1\cdot \hs1)$ is monotone increasing.

To verify {\rm\ci \eiiipxxxip}, we first notice that
$$
\int_{-k}^k\chi _{\tsize_{F^\star}}=\psi +k=\int_{-k}^k\chi _{\tsize_F},
\ct
$$ and hence
$$
\multline
\int_{-k}^ka_{\epsilon _0}(x,t)[\chi _{\tsize_{F^\star}}-\chi
_{\tsize_F}]=\int_{-k}^k[a_{\epsilon _0}(x,t)-a_{\epsilon _0}(x,\psi )][\chi
_{\tsize_{F^\star}}-\chi _{\tsize_F}]\\ =\int_{-k}^\psi [a_{\epsilon
_0}(x,t)-a_{\epsilon_0}(x,\psi )][1-\chi _{\tsize_F}]+\int_{\psi }^k[a_{\epsilon
_0}(x,t)-a_{\epsilon _0}(x,\psi )][0-\chi _{\tsize_F}].
\endmultline\tag{\nh.\en}
$$ But both integrals are non-positive due to the monotonicity of
$a_{\epsilon_0}(x,\hs1\cdot \hs1)$, and  {\rm\ci \tiiipvp} is proved.
\enddemo

The function $a_{\epsilon _0}$ is locally bounded in $\widetilde N$---in fact we
could modify it so that it would be globally bounded---, from which we
immediately infer that the boundary of any local minimizer $E$ of $\Cal F$ is {\it
almost minimal\/}, i.e. for any $\widetilde \Omega \Su \widetilde N$ there exists
$R>0$ and a constant $c$\hsc1, such that
$$
\int_{B_\rho (x)}|D\chi _{\tsize_E}|\le \int_{B_\rho (x)}|D\chi _{\tsize_F}|+c\rho
^{n+1}
\ct
$$
\we{3}{34}{0}%
for any $x\in \widetilde \Omega $, any $0<\rho <R$\hsc1, and any
$F$ with $F\sd E\Su B_\rho (x)$.

This definition is a special case of a more general one, where the second term on
the right-hand side of {\rm\ci \eiiipxxxivp} is supposed to grow with exponent
$(n+2\alpha )$, $0<\alpha <1$, cf. {\rm\Cite{\tii}}.

We note, that almost minimal boundaries in $\widetilde N$---or any other
$(n+1)$-dimensional Riemannian space---are also almost minimal in $\bold
R^{n+1}$, hence the regularity results proved in Euclidean space apply, i.e. the
reduced boundary of an almost minimal hypersurface is of class $C^{1,\alpha }$,
thus, in our case of class $C^{1,\frac{1}{2}}$, and the singular set is empty if $n\le
6$.
\head\hn. Proof of the main theorem
\endhead\equationno=0\proclaimno=0 The $C^{1,\alpha }$-estimates for almost
minimal boundaries yield uniform a priori estimates in the case of a sequence of
almost minimal boundaries satisfying the condition {\rm\ci \eiiipxxxivp} or its
more general variant uniformly. Moreover, assuming that {\rm\ci \eiiipxxxivp}
holds uniformly for a sequence of Caccioppoli sets $E_k\su \widetilde N$ which
converge locally to some limit set $E$, then, for any convergent sequence $x_k\in
\partial E_k$ with $x=\lim x_k$ we have $x\in \partial E$\hsc1; if in addition $x\in
\partial^\star E$ (the reduced boundary), then there exists $k_0$, such that
$x_k\in \partial^\star E_k$ for any $k\ge k_0$ and the unit normals at $x_k$
converge to the unit normal at $x$, cf. {\rm\Cite{\tii, Theorem 1}.}

In view of our assumption $n\le 6$, there are no singular points, i.e. $\partial^\star
E=\partial E$, and we conclude that the subgraphs $E_k=\sub \varphi _k$, where
$\varphi _k=\log u_k$, converge to $E=\sub \varphi $, $\varphi =\log u$\hsc1;
$\partial E$ is almost minimal, is of class $C^{1,\frac{1}{2}}$ and the mean
curvature of $M=\graph u$ in $N$ is equal to $f$, cf. Theorem 4.2 below. Hence,
$M$ and $\partial E$ are of class $C^{2,\alpha }$ for any $0<\alpha <1$. We
emphasize that only $M$ is smooth and not necessarily
$u$.

To complete the proof of {\rm\ci \tpiiip}, we have to show that $M$ is
homeomorphic to $S^n$ and that the mean curvature of $M$ is equal to $f$. For the
verification  of the spherical type of $M$ we observe that each
$M_k=\graph u_k$ is homeomorphic to
$S^n$ and that we have
\proclaim{Proposition \cp}For large $k$ the hypersurfaces $\partial E_k$ are
graphs over $\pa E$.
\endproclaim
\demo{Proof}Each point $x\in \pa E$ is the limit of a sequence of points $x_k\in
\pa{E_k}$, and the corresponding unit normals
$\nu _k$ converge uniformly to $\nu $. $\pa E$ is therefore oriented and, by
construction, all $\pa{E_k}$ lie on one side of
$\pa E$. Let $d$ be the signed distance function of $\pa E$; it is of class
$C^{2,\alpha }$ in a small tubular neighbourhood $\Cal U$ of $\pa E$. Then, only
finitely many $\pa{E_k}$'s are not completely contained in $\Cal U$. Fix $k$ such
that $\pa{E_k}\su\Cal U$, then for any $y\in\pa {E_k}$ there is exactly one $x\in\pa
E$ such that
$$
\text{dist}(y,x)=d(y)
\ct
$$
\we{4}{1}{0}%
 We claim furthermore, that, if $\Cal U$ is chosen small enough, any normal
geodesic starting at an arbitrary point $x\in\pa E$---and pointing in the right
direction---, intersects $\pa{E_k}$ in exactly one point, which together with
{\rm\ci \eivpip} yields that $\pa {E_k}$ is a graph over $\pa E$.

To verify that claim, let us consider normal Gaussian coordinates $(x^\alpha )$
relative to $\pa E$, where $x^0$ represents the coordinate axis normal to $\pa E$.
We also suppose that the unit normal $\nu $ of $\pa E$ has coordinates $\nu
=(1,0,\dots,0)$. Since the unit normals $\nu _k^{}$ of $\pa {E_k}$ converge
uniformly to $\nu $ we conclude, that $\nu _k^0$ is as close to
$1$ as we wish for large $k$, but then $\pa{E_k}$ is at least locally a graph over
$\{x^0=0\}$,  e.g., $\pa{E_k}=\graph \eta $ (locally), cf. {\rm\Cite{\giu, Proposition
4.9}}. But then
$$
\eta (x)=d(y).
\ct
$$ 
\enddemo

Finally, let us verify that the mean curvature in $N$ of the corresponding limiting
hypersurface $M$ is equal to $f$ which is not so obvious.

\proclaim{Theorem \cp}The mean curvature of $M$ is equal to $f$.
\endproclaim
\wt{4}{2}{0}%
To prove the theorem we need the following lemmata.

\proclaim{Lemma \cp}Let $K\su S^n$ be compact with $H^{n-1}(K)<\infty$, then
$$
H^n(\pa E\cap (K\times \bold R))=0.
\ct
$$
\endproclaim
\wl{4}{3}{0}%
\demo{Proof} We consider $\partial E$ and $\partial E_k$ as submanifolds of
$\widetilde N$. For a subset $U\su S^n$ we define $\widehat   U=U\ti \bold R$,
where we use polar coordinates. We also denote by $\mu$  resp. $\mu_ k$ the
measures
$|D\chi _E|$ resp. $|D\chi _{E_{_k}}|$, where $E$ resp. $E_k$ are the
subgraphs defined above.

Now, let $\Omega =S^n\smallsetminus K$; then $\Omega $ is a Caccioppoli set,
since $H^{n-1}(\partial \Omega )$ is finite. For $\epsilon >0$, let $\xi =(\xi ^\alpha
)\in H^{1,p}(\widetilde N)$ be a vectorfield such that

$$
|\xi |\le 1\qquad \text{ and}\qquad |\xi -\nu |_{_{\partial E}}\le \epsilon .
\ct
$$

Here, we use the fact that we already know $\partial E$ and $\partial E_k$ to be
uniformly of class $H^{2,p}$ for any finite $p$, since their mean curvatures are
uniformly bounded, and we have established a prioiri estimates in $C^{1,1/2}$.

Let $Q=S^n\ti [t_1,t_2]$ such that $\partial E_k \su Q$ for all $k$. Then we deduce

$$
\aligned
0&=\int_Q D_\alpha [\xi ^\alpha (\chi _{_{\ssize E_{k}}}-\chi _{_{\ssize E}})\chi
_{_{\widehat\Omega }}] \\
&=\int_Q D_\alpha \xi ^\alpha (\chi _{_{\ssize
E_{k}}}-\chi _{_{\ssize E}})\chi _{_{\widehat\Omega }} 
+ \int_Q \xi ^\alpha D _\alpha
\chi _{_{\ssize E_{k}}}\chi _{_{\widehat\Omega }}\\
&\phantom{+}-\int_Q \xi ^\alpha D _\alpha
\chi _{_{\ssize E}}\chi _{_{\widehat\Omega }}
+\int_Q  \xi ^\alpha (\chi _{_{\ssize E_{k}}}-\chi _{_{\ssize
E}})D_\alpha\chi _{_{\widehat\Omega }}\\
&\equiv I_1+I_2+I_3+I_4
\endaligned
\ct
$$
$I_1$ tends to zero if $k$ goes to infinity, since $\chi _{{\vphantom{}}_{E_k}}$
converges pointwise to $\chi _{_E}$. The same argument applies to $I_4$, while
$I_3$ is estimated by 
$$
|I_3|\le \mu (\widehat \Omega )
\ct
$$

Thus, we conclude
$$
\aligned
\varlimsup \mu _k(\widehat \Omega )&=\varlimsup \int_Q \chi
_{{\vphantom{}}_{\widehat \Omega }}|D\chi _{{\vphantom{}}_{E_k}}|\\
&\le \mu (\widehat \Omega )+c\,|\xi -\nu |_{{\vphantom{}}_M}\le \mu
(\widehat \Omega ) +c\,\epsilon ,
\endaligned
\ct
$$
because of the uniform convergence of $\nu_ k$ to $\nu $.

On the other hand, we have
$$
\mu (\widehat \Omega )\le \varliminf \mu_ k(\widehat \Omega ),
\ct
$$
and hence
$$
\mu (\widehat \Omega )=\lim \mu _k(\widehat \Omega ).
\ct
$$

Finally, we observe that $\mu _k(\widehat K)=0$ and
$$
\mu (\widetilde  N)=\lim \mu _k(\widetilde  N)
\ct
$$
and deduce the desired result.
\enddemo
In the next lemma we consider a  Caccioppoli set $E\su \widetilde N$
which is the subgraph of a function $\varphi \in BV(S^n)$.

\proclaim{Lemma \cp}Let $\varphi \in BV(S^n)$, $E=\sub \varphi $, and $(x,t)$ be
an interior point of $E$. Then, the line $\{(x,\tau)\,:\,-\infty<\tau<t\}$ does not
intersect the measure-theoretical boundary of $E$, i.e. the set of all points
$z$ such that
$$
0<|E\cap B_\rho (z)|<|B_\rho (z)|\qquad \forall \; 0<\rho <\rho (z).
\ct
$$
\endproclaim
\wl{4}{4}{0}%
\Cmedskip
Here, $B_\rho (z)$ is the geodesic ball of radius $\rho $ and center $z$, and
$|B_\rho (z )| $ its volume. We note that due to the metric product structure of
$\widetilde N$
$$
|B_\rho (x,t)|=|B_\rho (x,\tau ) |\quad \text{and}\quad \chi _{{\vphantom{}}_{B_\rho
(x,t) }}(\cdot ,\cdot +\tau )=\chi _{{\vphantom{}}_{B_\rho (x,t-\tau )
}}(\cdot ,\cdot )
\ct
$$
for arbitrary $t$ and $\tau $.
\demo {Proof of the Lemma}We denote the measure-theoretical boundary of $E$
by $\partial E$, since in the case we have in mind, the measure-theoretical and the
topological boundary coincide.

First, let us observe that the partial derivative of $\chi _{{\vphantom{}}_E}$ with
respect to $- \frac {\partial }{\partial x^0}$, $-D_0\chi _{{\vphantom{}}_E}$, or
more precisely,
$-\langle \frac {\partial }{\partial x^0}, D\chi _{{\vphantom{}}_E}\rangle$, is a
non-negative measure. For let $\eta \in C_c^\infty(\widetilde N)$, then
$$
\int_{S^n}\eta (x,\varphi (x))=\int_{S^n}\int_{-\infty}^{\varphi (x)}D_0\eta
(x,t)=-\int_{\widetilde N}\eta D_0\chi _{{\vphantom{}}_E}.
\ct
$$

Secondly, let $(x,\tau )\in \partial E$ with $\tau <t$. Then, we claim
$$
|E\cap B_\rho (x,\tau )|-|E\cap B_\rho (x,t)|=-\int_0^{t-\tau }ds\int_{B_\rho
(x,t-s)}D_0\chi _{{\vphantom{}}_E}.
\ct
$$
\we{4}{14}{0}%
The proof of this relation is exactly the same as that of its Euclidean counterpart,
cf. {\rm\Cite{\giu , Lemma 4.5}}.

Now, the right-hand side of {\rm\ci \eivpxivp} is non-negative, while the
left-hand side is strictly negative for small $\rho $, a contradiction.
\enddemo
We return to our original meaning of $\varphi $ and define
\proclaim{Definition \cp}Let $[\varphi ](x)$ be the jump of $\varphi $ at
$x$, i.e.
$$
[\varphi ](x)=\varlimsup_{y\rightarrow x}\varphi (y)-\varliminf_{y\rightarrow
x}\varphi (y)\equiv \varphi^+(x)-\varphi^-(x).
\ct
$$
\endproclaim
We have of course $\varphi ^+(x)=\varphi (x)$ since $\varphi $ is u.s.c.

An immediate corollary of {\rm\ci \livpivp} is

\proclaim{Lemma \cp}Let $[\varphi ](x)>0$, then the whole line segment $x\ti
[\varphi ^-(x),\varphi^+(x)]$ belongs to $\partial E$.
\endproclaim

The proof is straightforward since $(x,\varphi ^-(x))\in \partial E$.

\proclaim{Lemma \cp}Let $\tau >0$ and

$$
\Lambda _\tau =\{x\in S^n\,:\,[\varphi ](x)\ge \tau \}.
\ct
$$
Then, $\Lambda _\tau $ is compact and $H^{n-1}(\Lambda _\tau )\le c\,\tau ^{-1}$.
\endproclaim
\wl{4}{7}{0}%
\demo{Proof}We use a Besicovitch type covering argument. Let $0<\delta <\tau $,
then there exists a sequence of pairwise disjoint balls $B_{\rho _i}(x_i)$ in
$\bold R^{n+1}$, with centers $x_i\in \Lambda _\tau $ and radii $\rho
_i<\delta /3$, such that the balls $B_{3\rho _i}(x_i)$ cover $\Lambda _\tau $, see
e.g. {\rm\Cite{\giu, Lemma 2.2}}. We also choose $\delta $ small enough, such that
the volume of a geodesic ball in $\widetilde N$ of radius $\rho <\delta $ and center
in a compact set is uniformly bounded from below and above by a multiple of
$\rho ^{n+1}$. Consider the pairwise disjoint cylinders

$$
Q_{\rho _i}(x_i)=(B_{\rho _i}(x_i)\cap S^n)\ti \bold R.
\ct
$$

Then, $Q_{\rho _i}(x_i)\cap \partial E$ contains the line segment $x_i\ti [\varphi
^-(x_i),\varphi^+ (x_i)]$, and we can find $N_i$ disjoint geodesic balls
$B_{\rho_i/2}(y_{i,m}),\,1\le m\le  N_i$, with centers
$$
y_{i,m}\in x_i\ti [\varphi ^-(x_i),\varphi ^+(x_i)],
\ct
$$
where $N_i$ can be estimated by
$$
N_i\ge \frac{1}{4\rho _i}[\varphi ](x_i)\ge \frac{\tau }{4\rho _i}.
\ct
$$

Hence, we deduce
$$
\mu (Q_{\rho _i}(x_i))\ge \sum_{m=1}^{N_i}\mu (B_{\rho_i/2}(y_{i,m}))\ge
c\,\sum_{m=1}^{N_i}\rho _i^n\ge \frac{c}{4}\tau \rho _i^{n-1},
\ct
$$
where we use the well-known estimate for almost minimal boundaries
$$
\mu (B_{\rho_i/2}(y_{i,m}))\ge c\rho _i^n,
\ct
$$
if $y_{i,m}\in \partial E$, with a uniform positive constant $c$.

We infer further
$$
\sum_{i=1}^\infty\rho _i^{n-1}\le c\,\tau ^{-1}\sum_{i=1}^\infty\mu (Q_{\rho_
i}(x_i))\le c\,\tau ^{-1}\mu (\widetilde N),
\ct
$$
and conclude that the spherical $(n-1)$-dimensional measure of $\Lambda _\tau $
is bounded by a multiple of $\tau ^{-1}\mu (\widetilde N)$, but this is equivalent to
$$
H^{n-1}(\Lambda _\tau )\le c\,\tau ^{-1}\mu (\widetilde N)
\ct
$$
\we{4}{23}{0}%
with is a different constant.

$\Lambda _\tau $ is also closed, for let $x_m\in \Lambda _\tau $ be a sequence
converging to $x_0\in S^n$, then the line segments $x_m\ti [\varphi
^-(x_m),\varphi ^+(x_m)]$ in $\partial E$ converge to a line segment over $x_0$
of length at least $\tau $.
\enddemo

Combining {\rm\ci \livpiiip} and {\rm\ci \livpviip} we deduce
\proclaim{Lemma \cp}For each $\tau >0$ $H^{n-1}(\Lambda _\tau )=0$ and $\varphi
$ is $H^{n-1}$-a.e. continuous.
\endproclaim
\demo{Proof}We observe that for any Borel set $U\su S^n$
$$
\mu (\widehat U)=\int_U\sq\varphi .
\ct
$$
Hence, we have
$$
\int_{\Lambda _\tau }\sq\varphi =0,
\ct
$$
and for any $\epsilon >0$ there is an open set $\Omega $, $\Lambda _\tau \su
\Omega \su S^n$, such that
$$
\mu (\widehat \Omega )=\int_\Omega \sq\varphi <\epsilon .
\ct
$$

In the proof of {\rm\ci \livpviip} we can now choose the covering $B_{3\rho
_i}(x_i)$ such that
$$
B_{3\rho_i}(x_i)\cap S^n\su \Omega ,
\ct
$$
and instead of the estimate {\rm\ci \eivpxxiiip} we obtain
$$
H^{n-1}(\Lambda _\tau )\le c\,\tau ^{-1}\mu (\widehat \Omega) \le c\,\tau
^{-1}\epsilon ,
\ct
$$
i.e. $H^{n-1}(\Lambda _\tau )=0$.

The set where $\varphi $ is discontinuous is given by
$$
\bigcup_{k=1}^{\infty}\Lambda _{1/k},
\ct
$$
which is an $H^{n-1}$ null set.
\enddemo

We are now able to prove {\rm\ci \tivpiip}. The proof will be achieved, if we can
show
$$
\lim \int_{\partial E_k}(e^{\varphi _{k-1}}-e^{\varphi _k})=0.
\ct
$$

For large $k$ we can write $\partial E_k$ as a graph over $\partial E$
$$
\partial E_k=\{(\xi ,\eta _k):\,\,\xi \in \partial E\},
\ct
$$
where the $\eta _k$'s are uniformly of class $C^1$ and converge to $\eta =0$,
which corresponds to $\partial E$ in this setting. An integral of the form
$$
\int_{\partial E_k}f\; ,
\ct
$$
$f$ defined in $\widetilde N$, can then be expressed as
$$
\int _{\partial E}f(\xi ,\eta _k)\sq{\eta _k}\,b(\xi ,\eta _k),
\ct
$$
where the continuous volume forms $\sq{\eta _k}\,b(\xi ,\eta _k)$ converge to
$b(\xi ,0)$, the volume form for $\partial E$. This can be readily seen by
introducing normal Gaussian coordinates in a tubular neighbourhood of $\partial E$.

We extend $\varphi _{k-1}$ to $\widetilde N$ by the definition $\varphi
_{k-1}\circ x$, $x$ is the projection on $S^n$, and observe that $\varphi _k $ is
equal to $x^0_{| \partial E_k}$, where we still use polar coordinates $(x^i,x^0)$ in
$\widetilde N$.

Thus, we have
$$
\int_{\partial E_k}(e^{\varphi _{k-1}}-e^{\varphi _k})=\int _{\partial E}[e^{\varphi
_{k-1}\circ x(\xi ,\eta _k)}-e^{x^0(\xi, \eta _k)}]\sq{\eta _k}\,b(\xi ,\eta _k).
\ct
$$
Let $l\ge k_0$ be large, so that $\partial E_l$ is a graph over $\partial E$, and have
in mind that the sequence $\varphi _k$ is monotone falling. Then, we have for $k\ge
l-1$
$$
\multline
\int _{\partial E}[e^{\varphi
_{k-1}\circ x(\xi ,\eta _k)}-e^{x^0(\xi, \eta _k)}]\sq{\eta _k}\,b(\xi ,\eta _k)\\
\le \int _{\partial E}[e^{\varphi
_{l}\circ x(\xi ,\eta _k)}-e^{x^0(\xi, \eta _k)}]\sq{\eta _k}\,b(\xi ,\eta _k)
\endmultline
\tag{\nh.\en}
$$
and hence 
$$
\varlimsup \int_{\partial E_k}(e^{\varphi _{k-1}}-e^{\varphi _k})\le 
\int_{\partial E}[e^{\varphi _{l}\circ x(\xi ,0)}-e^{x^0(\xi ,0)}]\,b(\xi ,0),
\ct
$$
or, if we let $l$ tend to infinity,
$$
\varlimsup \int_{\partial E_k}(e^{\varphi _{k-1}}-e^{\varphi _k})\le 
\int_{\partial E}[e^{\varphi \circ x(\xi ,0)}-e^{x^0(\xi ,0)}]\,b(\xi ,0)
\ct
$$
But
$$
\varphi \circ x(\xi ,0)-x^0(\xi ,0)\le [\varphi ](x(\xi, 0))
\ct
$$
and in view of the preceding lemmata we know that $\mu $-a.e. the jump of
$\varphi $ is zero, i.e.
$$
\lim \int_{\partial E_k}(e^{\varphi _{k-1}}-e^{\varphi _k})=0.
\ct
$$

\immediate\closeout\lab

\Refs
\refstyle{C}
\widestnumber\key{99}
\ref 
\key \rn \wref{bk}
\by I. Bakelman and B. Kantor
\paper Existence of spherically homeomorphic hypersurfaces in Euclidean space
with prescribed mean curvature
\jour Geometry and Topology, Leningrad,\vol 1\yr 1974\pages 3--10 \miscnote 
\endref
\ref 
\key \rn \wref{bh}
\by H. Br\'ezis
\paper Equations et in\'equations non lin\'eaires dans les espaces vectoriel en
dualit\'e
\jour Ann. Inst. Fourier\vol 18\yr1968 \pages 115--175\miscnote 
\endref
\ref 
\key \rn \wref{DG}
\by E. De Giorgi
\paper Frontiere orientate di misura minima
\jour Sem. Mat. Scuola Norm. Sup. Pisa\vol \yr 1960--61\pages \miscnote 
\endref
\ref 
\key \rn \wref{DGcp}
\by E. De Giorgi, F. Colombini, and L.C. Piccinini
\book Frontiere orientate di misura minima e questioni collegate
\bookinfo \publ Scuola Norm Sup. Pisa\publaddr Pisa\yr 1972\pages 177\lang 
\endref
\ref 
\key \rn \wref{cgi}
\by C. Gerhardt
\paper Closed Weingarten hypersurfaces in space forms
\inbook Geometric Analysis and the Calculus of Variations
\ed J. Jost
\bookinfo\publ    International Press\publaddr Boston\yr 1996\pages 71--98\lang 
\endref
\ref 
\key \rn \wref{cgii}
\bysame
\paper Hypersurfaces of prescribed Weingarten curvature
\jour Math. Z.\vol \yr \pages \miscnote to appear
\endref
\ref
\key \rn \wref{cgiii}
\bysame
\paper Closed hypersurfaces of prescribed Weingarten curvature in Riemannian
manifolds
\jour J. Diff. Geom.\vol 43\yr1996 \pages 612--641
\endref
\ref 
\key \rn \wref{cgiv}
\bysame
\paper Hypersurfaces of prescribed mean curvature
\jour Math. Z.\vol 133 \yr1973 \pages 169--185\miscnote 
\endref
\ref 
\key \rn \wref{giu}
\by E. Giusti
\book Minimal surfaces and functions of bounded variation
\bookinfo Monographs in Math. Vol. 80\publ Birkh\"auser\publaddr
Boston-Basel-Stuttgart\yr 1984\pages 240\lang 
\endref
\ref 
\key \rn \wref{masi}
\by U. Massari
\paper Esistenza e regolarit\`a delle ipersuperfici di curvatura media assegnata
in $\bold R^n$
\jour Arch. Rat. Mech. Analysis\vol 55\yr 1974\pages 357--382\miscnote 
\endref
\ref 
\key \rn \wref{masii}
\bysame 
\paper Frontiere orientate di curvatura media assegnata in $L^p$
\jour Rend. Sem. Mat. Univ. Padova\vol53 \yr1975\pages 37--52\miscnote 
\endref
\ref 
\key \rn \wref{masM}
\by U. Massari and M. Miranda
\book Minimal surfaces of codimension one
\bookinfo Mathematics Studies Vol. 91\publ North-Holland\publaddr
Amsterdam-New York-Oxford\yr 1984\pages 242\lang 
\endref
\ref 
\key \rn \wref{masP}
\by U. Massari and L. Pepe
\paper Successione convergenti di ipersuperfici di curvatura media assegnata
\jour Rend. Sem. Mat. Univ. Padova\vol 53\yr1975 \pages53--68 \miscnote 
\endref
\ref 
\key \rn \wref{mm}
\by M Miranda
\paper Superfice cartesiani generalizzate ed insiemi di perimetri localmente
finito sui prodotti cartesiani
\jour Ann. Scuola Norm. Sup. Pisa, Serie III\vol 18\yr 1964\pages
515--542\miscnote 
\endref
\ref 
\key \rn \wref{tii}
\by I. Tamanini
\paper Boundaries of Caccioppoli sets with H\"older-continuous normal vector
\jour J. Reine Angew. Math.\vol 334\yr 1982\pages 27--39\miscnote 
\endref
\ref 
\key \rn \wref{tiii}
\bysame 
\paper Regularity results for almost minimal oriented hypersurfaces in $\bold
R^N$
\jour Quaderni Dipartimento Mat. Univ. Lecce\vol 1\yr 1984\pages
1--92\miscnote 
\endref
\ref 
\key \rn \wref{tw}
\by A.\,E. Treibergs and S.\,W. Wei
\paper Embedded hypersurfaces with prescribed mean curvature
\jour J. Diff. Geom.\vol 18\yr 1983\pages 513--521\miscnote 
\endref
\endRefs
\enddocument